\documentclass[11pt,a4paper]{article}

\usepackage{graphicx, amssymb, amsmath, amsthm, amsfonts, amsbsy, color, textcomp}
\usepackage{verbatim}
 
\DeclareMathOperator{\CT}{CT}

\DeclareMathOperator{\res}{\text{Res}}

\newtheorem{theorem}{Theorem}
\newtheorem{lemma}{Lemma}[section]

\newtheorem{definition}{Definition}

\newcommand{\Abar}{\overline{A}}
\newcommand{\Bbar}{\overline{B}}
\newcommand{\Cbar}{\overline{C}}

\newcommand{\bm}{\left ( \begin{matrix}}
\newcommand{\fm}{\end{matrix} \right ) }

\newcommand{\la}{\lambda}
\newcommand{\ka}{\kappa}
\newcommand{\w}{\omega}

\newcommand{\p}{\rho}

\newcommand{\Ao}{\overline{A}\,}
\newcommand{\Bo}{\overline{B}\,}
\newcommand{\Co}{\overline{C}\,}

\newcommand{\integer}{\mathbb{Z}}

\newcommand{\bd}{\hat{b}}		
\newcommand{\lad}{\hat{\lambda}}	
\newcommand{\defin}[1]{\textbf{#1}}		

\newcommand{\rP}{L}		
\newcommand{\recip}[1]{R}		

\newcommand{\cS}{\mathcal{S}}

\newcommand{\df}{:=}

\begin{document}
%
%

\title{Chebyshev Type Lattice path weight polynomials by a constant term method}
\author{R.~Brak${}^\dag$\ and J.~Osborn${}^{\dag\dag}$\\
${}^\dag$Department of Mathematics and Statistics,\\
The University of Melbourne,\\
Parkville, Victoria 3010, Australia.\\
${}^{\dag\dag}$Centre for Mathematics and its Applications, \\
Mathematical Sciences Institute,\\
Australian National University, \\ 
Canberra, ACT 0200, Australia.
}

\maketitle

\begin{abstract}
We prove a constant term theorem which is useful for finding weight
polynomials for Ballot/Motzkin paths in a strip with a fixed number of
arbitrary `decorated' weights as well as an arbitrary `background'
weight.
Our CT theorem, like Viennot's lattice path theorem from which it is
derived primarily by a change of variable lemma,  is expressed in
terms of orthogonal polynomials which in our applications of interest
often turn out to be non-classical.
Hence we also present an efficient method for finding explicit closed
form polynomial expressions for these non-classical orthogonal
polynomials.
Our method for finding the closed form polynomial expressions relies
on simple combinatorial manipulations of  Viennot's diagrammatic
representation for orthogonal polynomials.
In the course of the paper we also provide a new proof of Viennot's
original orthogonal polynomial lattice path theorem.  The new proof is
of interest because it uses diagonalization of the transfer matrix,
but gets around difficulties that have arisen in past attempts to use
this approach.  In particular we show how to sum over a set of
implicitly defined zeros of a given orthogonal polynomial, either by
using properties of residues or by using partial fractions.
We conclude by applying the method to two lattice path problems
important in the study of  polymer physics as models of steric
stabilization and sensitized flocculation.

\end{abstract}


\section{Introduction and Definitions} 

As is well known to mathematical physicists the form of the solution to a problem is often more important than its existence. Such is the case in this paper. Determining the generating  function of   Motzkin path weight polynomials in a strip was solved by a theorem  due to Viennot \cite{Viennot1985ah} (see Theorem \ref{thm:viennot} below -- hereafter referred to as Viennot's Theorem). In the applications discussed below what is required are the weight polynomials themselves. Whilst these can be written as a Cauchy integral of the generating function this form of the solution is of little direct use for our applications. To this end we have derived a  related form  of the generating function given by Theorem \ref{thm:CT} (hereafter referred to as the Constant Term, or CT, theorem).  The CT theorem  is certainly well suited to extracting the weight polynomials for Chebyshev type problems (section \ref{sec:applications}) and as a starting point for their asymptotic analysis \cite{owczarek:2008qe}.  

The CT theorem is  proved, as we shall show  below, by starting with Viennot's Theorem  and using a `change of variable' lemma (see Lemma~\ref{lem:changeVars}). We will  also provide a new   proof of Viennot's Theorem that is based on diagonalizing the associated Motzkin path transfer matrix. The latter proof is included as it rather naturally leads to the CT Theorem. It also has some additional interest as it has several combinatorial connections \cite{brak:2006kx}. For example, a combinatorial interpretation of what has previously appeared only as a change of variable to eliminate a square root in Chebyshev polynomials turns out to be the generating function of binomial paths. Another connection is a combinatorial interpretation of the Bethe Ansatz \cite{bethe31} as determining the signed set of path involutions,  as  for example, in the involution of Gessel and Viennot \cite{Gessel:gr} in the many path extension \cite{brak:1998un}.

Two classes of  applications for which  the CT theorem   is certainly suited are the Asymmetric Simple Exclusion Process (ASEP) and directed models of polymers interacting with surfaces.  For the ASEP the problem of computing the stationary state,  and hence in finding phase diagrams for associated simple traffic models, can be cast as a lattice path problem \cite{blythe:2007yq,blythe:2000ru, blythe:ai,  brak:2006fj, Brak2004vf,Corteel05ve, derrida97, evans:1999cr}. For the  the ASEP model  the path problem required is actually a half plane model with two weights associated with the lower wall (the upper wall is sent to infinity to obtain the half plane.) \cite{Brak2004vf}.  
In chemistry the lattice  paths are used to model polymers in solution \cite{gennes79} - for instance in the analysis of steric stabilization and  sensitized flocculation  \cite{brak:2007lr,brak:2005fa}. 

In the application section of this paper we find   explicit expressions for the partition function -- or weight polynomial -- of the DiMazio and Rubin polymer model \cite{dimarzio:1971ev}. This model was first posed in 1971  and has    Boltzmann weights   associated with upper and lower walls of a strip containing a path. The wall weights model the interaction of the polymer with the surface.  The solution given in this paper is an improvement on that published in \cite{brak:2006kx}. Previous results on the DiMazio and Rubin model have only dealt with special cases of weight values, for example, the case  where a  relationship exists between the Boltzmann weights \cite{brak:1999ri}.   We also present a natural generalisation of the DiMazio and Rubin weighting which may have application to models of polymers interacting with colloids \cite{owczarek:2007yq}, where the interaction strength depends on the proximity to the colloid. 

In order to use the CT Theorem both the ASEP and polymer models require explicit expressions for  `perturbed' Chebyshev orthogonal  polynomials. Their computation is addressed by  our third theorem, Theorem~\ref{thm:ViennotPaving}


\subsection{Definitions and Viennot's Theorem} 
\label{sec:definitions_and_viennots_theorem}

\begin{figure}[htbp]
\begin{center}
\includegraphics[width=12cm]{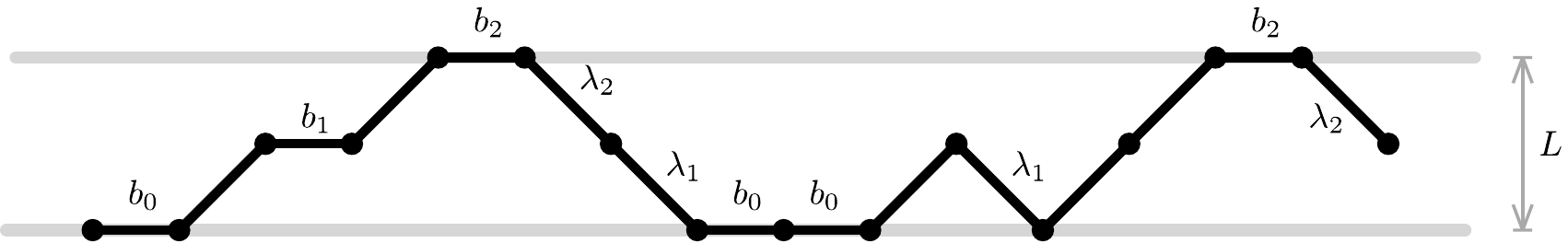}
\caption{An example of a lattice path of length 15, in a strip of height $L=2$, with weight $b_0^3b_1b_2^2\lambda_1^2\lambda_2^2$ and starting at $y'=0$ and ending at $y=1$.}
\label{fig:CTschemaRotated}
\end{center}
\end{figure}

Consider length $t$ lattice paths, $p= v_0v_1v_2 ... v_t$, in a height $L$ strip with vertices $v_i\in\cS=\mathbb{Z}_{\ge 0} \times \{0,1,...,L\}$, such that the  edges $e_i \df v_{i-1} v_i$ satisfy $v_i - v_{i-1}\in\{(1,1), (1,0),(1,-1) \}$.  An edge $e_i$  is called \defin{up} if $v_i - v_{i-1}=(1,1)$ , \defin{across} if $v_i - v_{i-1}=(1,0)$ and \defin{down} if $v_i - v_{i-1}=(1,-1)$. A vertex $v_i=(x,y)$ has \textbf{height} $y$.  Weight the edges according to the height  of their left vertex using
\begin{equation} \label{eq:generalpathWeights}
w(e_i) =
\begin{cases}
	1 & \text{if $e_i$ is an \emph{up} edge} \\
	b_{k} & \text{if $e_i$ is an \emph{across} edge with $v_{i-1}=(i-1, k)$}\\
	\la_k & \text{if $e_i$ is a \emph{down} edge with $v_{i-1}=(i-1, k)$}.
\end{cases}
\end{equation}
The weight of the path, $w(p)$, is defined to be the product of the weights of the edges, ie.\ for path $p= v_0   v_1 ... v_t$,
\begin{equation}
w(p) = \prod_{i=1}^t w(e_i).
\end{equation}
Such weighted paths are then enumerated according to their length  with \defin{weight polynomial}  defined by
\begin{equation}\label{eq:weightpoly1}
Z_t(y', y;L) \df  \sum_{p} w(p) 
\end{equation}
where the sum is over all paths of length $t$, confined within the strip of height $L$,  $y'$ is the  height of the initial vertex of the path and $y$ is the height  of the  final vertex of the path. An example of such a path is shown in Figure~\ref{fig:CTschemaRotated}.
These paths are weighted and confined elaboration's of Dyck paths, Ballot paths and Motzkin paths -- for enumeration results on these classic paths see, for example, \cite{goulden:1983yg}.  Once non-constant weights are added, many classical techniques do not (obviously) apply.

This work focuses upon solving the enumeration problem for the types of weighting in which a small number of weights take on distinguished values called `decorations', and the rest of the edges have constant `background weights' of just one of two kinds, accordingly as the step is an across step or a down step.

We introduce notation to describe the positions of the decorated  edges:
let $\mathcal{D}_{b} \subseteq \{0,...,L\}$ and $\mathcal{D}_{\la} \subseteq \{1,...,L\}$ be sets of integers called respectively \defin{decorated across-step heights} and \defin{decorated down-step heights}.  Then paths are weighted as in Equation~(\ref{eq:generalpathWeights}), with
\begin{subequations}\label{eq:lailall} \begin{align} 
	b_i &= \begin{cases}
			b 				& 	\text{if $i \notin \mathcal{D}_{b}$}\\
		  b + \bd_i &  \text{if $i \in\mathcal{D}_{b}$}
	\end{cases}
	 \\
	\la_i &= \begin{cases}
			\la 				& 	\text{if $i \notin \mathcal{D}_{\la}$}\\
		  \la + \lad_i &  \text{if $i \in\mathcal{D}_{\la}$}
	\end{cases} \label{eq:laila}	
\end{align}\end{subequations}
where  $b$ and $\la$ will be called \defin{background weights} and $\bd_i$ and $\lad_i$ will be called \defin{decorations}.

The generating function for the weight polynomials is given in terms of orthogonal polynomials by a theorem due to Viennot\cite{Viennot1985ah,viennot:1984fr}  -- see also \cite{Flajolet1980rr,goulden:1986kx}
\begin{theorem}[\cite{Viennot1985ah}]\label{thm:viennot}
	The generating function of the weight polynomial \eqref{eq:weightpoly1}  is given by
	\begin{equation}\label{eq:vienratio}
		M_L(y', y;x) \df \sum_{t\ge0} Z_t(y', y;L)\,x^t =x^{Y-Y'}\,\frac{\recip{P}_{Y'}(x)\,h_{y',y}\, \recip{P}^{(y+1)}_{L-Y}(x)}{\recip{P}_{L+1}(x)}, 
	\end{equation}
 where $h_{y',y}=  \prod_{y < l \le y'} \la_l$ if $y'> y$ and $h_{y',y}=1$ if $y'\le y$, $Y' =  \min \{ y', y\}$ and $ Y =  \max\{ y', y\}$. The polynomials  
  $\recip{P}_k(x)$  are the \textbf{reciprocal} polynomials and    $\recip{P}^{(j)}_k(x)$ the \textbf{shifted reciprocal}  polynomials  defined by
	\begin{equation}\label{eq:recippolys}
		\recip{P}_k(x)\df x^k P_k(1/x)\qquad\text{and}\qquad 
		\recip{P}^{(j)}_k(x)\df \biggl.\recip{P}_k(x)\biggr|_{\stackrel{b_i\to b_{i+j}}{\la_i\to \la_{i+j}}} 
	\end{equation}
	where the orthogonal polynomials $P_k(x)$ satisfy the standard three term recurrence \cite{chihara:1978lr,szego:1975ep},
	\begin{equation}\label{eq:threeTerm1}
		P_{k+1}(x)=(x-b_k)P_{k}(x)-\la_k P_{k-1}(x),\qquad P_0(x)=1,\quad P_1(x)=x-b_0.
	\end{equation}
\end{theorem}
This theorem may be proved in several ways; as the ratio of determinants -- see  \cite{stanley:1997vw} section 4.7.2,  by continued fractions \cite{Flajolet1980rr, goulden:1986kx} or by heaps of monomers and dimers \cite{viennot:1986fk}. In Section~\ref{sec_tmatrixPrrof} we provide a new proof that uses   diagonalization of the transfer matrix of the Motzkin paths.   


\section{A Constant Term Theorem} 
\label{sec_main_results}

Our  main result is stated as a constant term of a particular Laurent expansion.  Since the constant term method studied in this paper depends strongly on the choice of Laurent expansion, we briefly recall a few simple facts about Laurent expansions. Since we only consider rational functions  we restrict our discussion to them. A Laurent expansion of a rational function about a point   $z=z_i$ is of the form $\sum_{n\ge n_0} a_n(z-z_i)^n$. The coefficients $a_n$ depends  on the chosen point $z_i$ and the annulus of convergence. 
Furthermore, the nature of $n_0$ generally depends on three factors: i) $n_0\ge 0$ (ie.\ the series is a Taylor series) if the annulus contains no singular points, ii) $n_0<0$ is finite if the inner circle of the annulus contains only a non-essential singularity at $z_i$ and iii)  $n_0=-\infty$ if the inner circle contains at least one other singularity at $z\ne z_i$ or an essential singularity. 

In this paper we only need case ii), with $z_i=0$ and the series  convergent in the  annulus closest to the origin.  Thus the constant term is defined as follows.
  
\begin{definition} \label{def:CT} Let $f(z)$ be a complex valued function with Laurent expansion of the form
\begin{equation}
f(z) = \sum_{n=n_0}^{\infty} a_n z^n
\label{eq_laurentseries}
\end{equation}
with $n_0\in\integer$.
Then the \defin{constant term in $z$} of $f(z)$ is
\begin{equation}
\CT_z [f(z)] = a_0.
\end{equation}
\end{definition}  
This is, of course, just the residue of $f(z)/z$ at $z=0$.  Note, the form of the Laurent expansion given in \eqref{eq_laurentseries} uniquely specifies that it corresponds to that Laurent expansion of $f(z)$ that converges in the innermost annulus   that is centred at the origin.

Our  main result gives the weight polynomial for  Motzkin paths in a strip as the constant term  
 of a rational function constructed from Laurent  polynomials. The Laurent polynomials we use, $\rP_k^{(j)}(\p)$, are defined in terms of the conventional (shifted) orthogonal  polynomials, $P_k^{(j)}(x)$, by the simple substitution
\begin{equation}\label{eq:PR}
\rP_k^{(j)}(\p) \df  P_k^{(j)}\bigl(x(\p)\bigr)
\end{equation}
with  
\begin{equation} \label{eq:mup}
x(\p) =  \p+ b + \la \p^{-1}.
\end{equation}
The orthogonal  polynomials $P_k(x)=P_k^{(0)}(x)$ satisfy the standard three term recurrence \eqref{eq:threeTerm1} which, for the shifted polynomials $P_k^{(j)}(x)$, becomes
\begin{align} \label{eq:threetermrec}
	P_{k}^{(j)}(x) &= (x - b_{k+j-1})P_{k-1}^{(j)}(x) - \la_{k+j-1} P_{k-2}^{(j)}(x),\quad k\ge 2\\
	P_1^{(j)}(x) &= x-b_j\notag\\
	P_0^{(j)}(x) &= 1\notag
\end{align}
and $\la_k\ne0$ $\forall k$. We   now state our principal theorem. 
\begin{theorem}[Constant Term]
\label{thm:CT} 
	Let $Z_t(y',y;L)$ be the weight polynomial for the set of Motzkin paths with initial height $y'$, final height $y$, confined in a strip of height $L$, and weighted as specified in Equations~\eqref{eq:lailall}.  Then
\begin{equation} \label{eq:CTrho}
	Z_t(y',y;L)= \CT_\p \left [ \left (\p+ b +  \frac{\la}{\p} \right)^t 
	\frac{ \rP_{Y'}(\p)\, h_{y',y}\, \rP_{L-Y}^{(Y+1)}(\p) }{ \rP_{L+1}(\p)} \left (\frac{\la}{\p}-\p\right )\right ], 
\end{equation}
with $Y' =  \min \{ y', y\}, Y =  \max\{ y', y\}$, 
\begin{equation}
	h_{y',y}=\begin{cases}
			\prod_{y < l \le y'} \la_l & \text{if $y'> y$}\\
			1 & \text{otherwise.}
	\end{cases}
	\label{eq_hfactor}
\end{equation}
and the Laurent polynomials $\rP_k^{(j)}$ given by \eqref{eq:PR}.
\end{theorem}

The form of this constant term expression should be carefully compared with that arising from Viennot's theorem when used in conjunction with the standard Cauchy constant term form for the weight polynomial (for $y'\le y$),
\begin{equation}\label{eq:classCT}
	Z_t(y',y;L)= \CT_x\left[\frac{1}{x^{t+1}}x^{y-y'}\,\frac{\recip{P}_{y'}(x)\,  \recip{P}^{(y+1)}_{L-y}(x)}{\recip{P}_{L+1}(x)}\right].
\end{equation}
In particular, \eqref{eq:CTrho} is \emph{not} obtained by simply substituting $1/x=\p+ b + \la \p^{-1}$ into \eqref{eq:classCT} as was done to define the Laurent polynomials $\rP_k^{(j)}$.

It is a useful exercise to  compare the difference in effort in  computing a \emph{general expression} for $Z_t(y',y;L)$ in the simplest possible case $y=y'=0$,  $b_k=1$, $\la_k=\la$  starting from \eqref{eq:CTrho} compared with starting from \eqref{eq:classCT}.



Assuming that a simple explicit expression is desired for the weight polynomial then the utility of the CT Theorem depends on and arises from three factors. The first problem is how to calculate the orthogonal polynomials. For various choices of the weights    $b_k$ and $\la_k$, the classical orthogonal polynomials are obtained and hence this problem has already been solved. However  for the applications mentioned earlier, which require `decorated' weights, the polynomials\footnote{They may  be thought of as `perturbed' Chebyshev polynomials} do not fall into any of the classical classes. Thus computing the polynomials becomes a problem in itself and is addressed by Theorem~\ref{thm:ViennotPaving}.

The second problem is more subtle and is concerned with how the polynomials are represented. Whilst $P_{k}(x)$ is by construction a  polynomial  this is not necessarily how it is first represented. For example, 
consider   a  Chebyshev type polynomial,
which satisfies the constant coefficient recurrence relation
\begin{equation} \label{eq:Sk1}
S_{k+1}(x)=x S_k(x)- S_{k-1}(x),\qquad S_1(x)=x,\qquad S_0(x)=1.
\end{equation}
This recurrence is easily solved by substituting the usual trial solution $S_k=\nu^k$ with $\nu$ a constant, leading immediately to a solution in  the form
\begin{equation} 
S_k(x)=\frac{(x+\sqrt{x^2-4})^{k+1}-(x-\sqrt{x^2-4})^{k+1}}{2^{k+1}\sqrt{x^2-4}}.
\label{eq_chebEx}	
\end{equation}
Whilst this is a polynomial in $x$, in this form it is not \emph{explicitly} a polynomial\footnote{The square roots can of course be  Taylor expanded to show explicitly it is a polynomial.} as it is written in terms of the branches of an algebraic function. The representation of the polynomials is important as it strongly influences the third problem, that of computing the constant term (or residue).  If the polynomials are explicitly polynomials (rather than, say, represented by algebraic functions) then the obvious way of computing the Laurent expansion, and hence residue, is via a geometric expansion of the denominator. Whilst in principal this can always be done for a rational function the simpler the denominator polynomials the simpler the weight polynomial expression -- in particular we would like as few summands as possible, preferably a number that does not depend on  $L$, the height of the strip.
  
The fact that this can be achieved for the applications studied here  shows the advantage of the CT  Theorem in this context  over, say, the Rogers Formula \cite{rogers:1907rc} (see Proposition 3A of \cite{Flajolet1980rr}), in which the weight polynomial is always expressed as an $L$-fold sum no matter how small the set of decorated weights.

\section{A Paving Theorem} 
\label{sec:paving_theorem}

Our second theorem    is  used to find explicit expressions for the orthogonal polynomials which are useful to our applications.  These are polynomials arising from problems where the number of decorated weights is fixed (ie.\ independent of $L$). 
Theorem~\ref{thm:ViennotPaving}, of which we make extensive use,
expresses the orthogonal polynomial of the decorated weight problem in terms of the orthogonal polynomials of the problem with no decorated weights (ie.\ Chebyshev type polynomials).

\begin{theorem}[Paving] \label{thm:ViennotPaving} 
\begin{enumerate}

	\item \label{it:oneCut} For each $c \in \{1,2,...,k-1\}$, we have an  `edge cutting' identity
\begin{subequations}
\begin{equation}
P_k^{(j)}(x)=P_c^{(j)}(x)P_{k-c}^{(j+c)}(x)
	-\la_{c+j}\, P_{c-1}^{(j)}(x)P_{k-c-1}^{(j+c+1)}(x) \label{eq:edgeIdentity}
\end{equation}	
and a `vertex-cutting' identity,
\begin{multline}
P_k^{(j)}(x)=(x-b_{c+j})\,P_c^{(j)}(x)P_{k-c-1}^{(j+c+1)}(x)\\
-\la_{j+c+1}\,P_c^{(j)}(x)P_{k-c-2}^{(j+c+2)}(x)
	-\la_{c+j}\, P_{c-1}^{(j)}(x)P_{k-c-1}^{(j+c+1)}(x), 
	\label{eq:vertexIdentity}
\end{multline}	
\end{subequations}
where $P_k^{(j)}(x)$ satisfies \eqref{eq:threetermrec}.
	\item \label{it:bound} 
Fix $j$ and $k$.  Let $|\mathcal{D}_b|$ and $|\mathcal{D}_\lambda|$ be the number of decorated `across' and `down' steps respectively whose indices are  strictly between $j-1$ and $j+k$.  Let $d=|\mathcal{D}_\lambda|+|\mathcal{D}_b|$.  Then
\begin{equation}\label{eq:pertexp}
P_k^{(j)}(x)=	\sum_{j=1}^{j_{\text{max}}} a_j \prod_{i=1}^{i_{\text{max}}} S_{k_{j,i}}(x),
\end{equation}
where 
\begin{equation} \label{eq:bounds}
1 \le j_{\text{max}} \le 2^{|\mathcal{D}_\lambda|}3^{|\mathcal{D}_b|}, \,\,\,\,\,\, 1 \le i_{\text{max}} \le d+1;
\end{equation}
and $k_{j,i}$ is a positive integer valued function; decorations are all contained in the coefficient $a_j$'s, and the $S_{k_{j,i}}(x)$'s are the background weight dependent Chebyshev  orthogonal polynomials satisfying
\begin{equation} \label{eq:Sk}
S_{k+1}(x)=(x-b)\, S_k(x)- \lambda\, S_{k-1}(x),\quad 
\end{equation}
with $S_1(x)=x-b$, $S_0(x)=1$ and   $\la\ne 0$.
\end{enumerate}
\end{theorem}

We do not give an explicit expression  for the $k_{j,i}$ as it is strongly dependent on the sets $\mathcal{D}_b$ and $\mathcal{D}_\lambda$. They are however simple to compute in any particular case, for example see \eqref{eq:dmPoly2new} in the application section \ref{sec:apps}. The significance of \eqref{eq:pertexp} is that it shows that the decorated polynomials $P_k^{(j)}$ can be explicitly expressed in terms of the undecorated (ie.\ Chebyshev) polynomials $S_{k}^{(m)}(x)$

The first part of the paving Theorem~\ref{thm:ViennotPaving} follows immediately from an `edge-cutting' and a `vertex-cutting' technique respectively, applied to Viennot's paving representation of orthogonal polynomials, which we describe in  Section \ref{sec:applications}.  This geometric way of visualising an entire recurrence in one picture is powerful; from it we see Part~\ref{it:oneCut} of the Theorem as a gestalt, so in practice do not need to remember the algebraic expressions but may work with  paving diagrams directly.  Part~\ref{it:bound} follows immediately by induction on Part~\ref{it:oneCut}.


%
 
\section{Proof of Viennot's theorem by Transfer Matrix Diagonalization}	
\label{sec_tmatrixPrrof}

In this section we state a new proof of Viennot's theorem.   This proof starts with the transfer matrix for the Motzkin path (see section 4.7 of \cite{stanley:1997vw} for an explanation of the transfer matrix method) and proceeds by diagonalizing the matrix. As is well known this requires summing an  expression over all the eigenvalues of the matrix. The eigenvalue sum is a sum over the zeros of a particular orthogonal polynomial. This sum can be done for the most general orthogonal polynomial even though the zeros are not  explicitly known.  We do this in two ways, the first uses two classical results (Lemma~\ref{lem:asdasd} and Lemma~\ref{lem:insideoutside}) and a paving polynomial identity. The essential idea is to replace the sum by a sum over residues and this residue sum can then be replaced by a single residue at infinity. The second proof uses partial fractions.  

The \defin{transfer matrix}, $T_L$, for paths in a strip of height $L$, is a square matrix of order $L+1$ such that the $(y', y)^{\text{th}}$ entry of the $t^{\text{th}}$ power of the matrix gives the weight polynomial for paths of length $t$, i.e. 
\begin{equation}
Z_t(y',y;L) = (T_L^t)_{y', y}.
\end{equation}
For Motzkin paths with weights \eqref{eq:lailall} the transfer matrix is the Jacobi matrix,
\begin{equation} \label{eq:transferMatrix}
T_L \df  
\left ( \begin{matrix}
b_0 & 1 & 0 &  \cdots &  & &  0\\
\la_1 & b_1 & 1 & 0 &  \cdots  & & 0\\
0 & \la_2 & b_2 & 1 & 0 & \cdots & 0\\
\vdots  & \ddots & \ddots & \ddots & \ddots & \ddots & \vdots\\
0 & \cdots & 0 & \la_{L-2} & b_{L-2} & 1 & 0\\
0 & & \cdots &  0 & \la_{L-1} & b_{L-1} & 1\\
0 & & & \cdots &  0 & \la_{L} & b_{L} \\
\end{matrix} \right ) .
\end{equation}

The standard  path length generating function  for such paths, with specified initial height $y'$ and final height  $y$, is given  in terms of powers of the transfer matrix as
\begin{equation} \label{eq:GxT}
			M_L(y', y;x) \df  \sum_{t \ge 0}\,  Z_t(y',y;L)\, x^t = \sum_{t \ge 0} (T_L^t)_{y',y}x^t  
\end{equation}
which is convergent for $|x|$ smaller than the reciprocal of the absolute value  of the largest eigenvalue of $T_L$.

The details of the proof are as follows. We evaluate $Z_t(y',y;L)=(T_L^t)_{y',y}$ by  diagonalization
\begin{equation} \label{eq:TVDU}
T_L^t = V D_L^t U
\end{equation}
with $D_L=\text{diag}(x_0, x_1, ..., x_L)$ a diagonal matrix of eigenvalues of $T_L$, and $V$ and $U$ respectively matrices of right eigenvectors as columns and left eigenvectors as rows normalized such that
\begin{equation}
UV=I,
\label{eq_completeness}
\end{equation}
where $I$ is the unit matrix.
One may check that this diagonalization is achieved by setting the $i^\text{th}$ column  of $V$ equal to the transpose of
\begin{equation}
\mathbf{v}(x_{i-1}) = \left ( P_0(x_{i-1}), P_1(x_{i-1}), ..., P_L(x_{i-1})\right ) 
\end{equation}
and the $i^\text{th}$ row  of $U$ equal to
\begin{align}
\mathbf{u}(x_{i-1})=&\frac{\la_1 ... \la_L}{P'_{L+1}(x_{i-1})P_L(x_{i-1})}\notag\\ 
&\times\left (\frac{P_0(x_{i-1})}{1}, \frac{P_1(x_{i-1})}{\la_1}, \frac{P_2(x_{i-1})}{\la_1 \la_2}, ..., \frac{P_L(x_{i-1})}{\la_1 ... \la_L} \right),
\end{align}
where the set of eigenvalues $\{x_i\}_{i=0}^L$ are determined by
\begin{equation}
P_{L+1}(x_i)=0.
\label{eq_zeros}
\end{equation}
with the orthogonal polynomial $P_{L+1}$ given by the three term recurrence \eqref{eq:threetermrec}.
Orthogonality of left with right eigenvectors of the Jacobi matrix \eqref{eq:transferMatrix} follows by using the Christoffel-Darboux theorem for orthogonal polynomials. Equation \eqref{eq_completeness} then follows as \eqref{eq_zeros} gives $L+1$ distinct zeros and hence $L+1$ distinct eigenvalues and hence $L+1$ linearly independent eigenvectors.

For simplicity in the following we only consider the case $y'\le y$ in which case the $h_{y',y}$ factor is one -- it is readily inserted for the case $y' > y$. Thus, multiplying out Equation~(\ref{eq:TVDU}) and extracting the $(y', y)^{\text{th}}$ entry, we have
\begin{equation}
(T_L^t)_{y',y} =(\la_{y+1} ... \la_L) \sum_{i=0}^L \frac{x_i^t P_{y'}(x_i)P_y(x_i) }{P'_{L+1}(x_i)P_L(x_i)}. 
\label{eq_quadortho}
\end{equation}
Note that $P'_{L+1}(x_i) \ne 0$ and $P_L(x_i)\ne 0$ by the Interlacing Theorem for orthogonal polynomials, so that all the terms in the sum are finite. 
Since Viennot's theorem does not have a product of polynomials in the denominator we need to simplify \eqref{eq_quadortho}, which is achieved by using the following lemma.

%
\begin{lemma} 
	\label{lem:nontrivialPavingPolyIdentity} 
	Let $x_i$ be a zero of  $P_{L+1}(x)$.  Then
\begin{equation}
		\la_{y+1}...\la_L\,	P_y(x_i)
			=  P_L(x_i) P^{(y+1)}_{L-y}(x_i).
\end{equation}
\end{lemma} 
%
This lemma follows directly from the edge-cutting identity  \eqref{eq:edgeIdentity} by choosing $k=L+1$, $j=c$ and $c=L-h$ together with the assumption that  $P_{L+1}(x_i)=0$, to obtain a family of identities parametrized by $h$; which are then iterated with $h=0,1,\dots,L-y-1$.

Applying Lemma~\ref{lem:nontrivialPavingPolyIdentity}  to \eqref{eq_quadortho} gives the following basic expression for the weight polynomial resulting from the transfer matrix
\begin{equation}\label{eq:evexpand}
(T_L^t)_{y',y} =\sum_{i=0}^L \frac{x_i^t P_{y'}(x_i) P^{(y+1)}_{L-y}(x_i) }{P'_{L+1}(x_i)}. 
\end{equation}
Note, this use of the transfer matrix leads \emph{first} to an expression for the weight polynomial,
to get to Viennot's Theorem we still need to generate on the path length and also simplify the sum over zeros. The former problem is trivial, the latter not. We do the sum over zeros in two ways, first by using a contour integral representation and secondly by using partial fractions.

\subsection*{Using a contour integral representation}

To eliminate the derivative in the denominator of \eqref{eq:evexpand} we use the following Lemma.
\begin{lemma}  
\label{lem:asdasd}
Let $P(z)$ and $Q(z)$ be  polynomials in a complex variable $z$ and  $P(z_i)=0$, then 
\begin{equation}
	m_i\,	Q(z_i)= \res \left [ Q(z) \frac{P'(z)}{P(z)}, \{z, z_i\}\right ].
		\label{eq_zeroToRez}
\end{equation}
where $m_i$ is the multiplicity of the root $z_i$, and $P'(z)$ is the derivative with respect to $z$.
\end{lemma} 
%
The following Lemma allows us to replace a residue sum with a residue at infinity.
\begin{lemma}  
\label{lem:insideoutside}
Let $P(z)$ and $Q(z)$ be  polynomials in a complex variable $z$, then
\begin{subequations}
\begin{align}
	\frac{1}{2\pi i }\int_\gamma \frac{P(z)}{Q(z)}dz	&=	\sum_{z_i \in A} \res \left [ \frac{P(z)}{Q(z)}, \{z, z_i \}\right ]\\
		&=	 \res \left [\frac{P(1/z)}{z^2Q(1/z)}, \{z, 0\}\right ],
\end{align}		
\label{eq_insideoutside}
\end{subequations}
where $\gamma$ is a simple closed anticlockwise-oriented contour enclosing all the  zeros of $Q(z)$ and $A$ is the set of zeros of $Q(z)$. 
\end{lemma} 
%

Note, \eqref{eq_insideoutside} simply states that the sum of all residues of a rational function, including that at infinity, is zero.  These lemmas are proved in most books on complex variables -- see \cite{boas:1987lq} for example.

Now use Lemma~\ref{lem:asdasd} to get rid of the derivative in the denominator ($m_i=1$ as all zeros of orthogonal polynomials are simple), to produce
\begin{equation} 
(T_L^t)_{y',y} = \sum_{\{x_i |P_{L+1}(x_i)=0 \}} \res \left [ \frac{x^t P_{y'}(x)P^{(y+1)}_{L-y}(x)}{P_{L+1}(x)} , \{x, x_i\} \right].  
\end{equation}
Applying Lemma~\ref{lem:insideoutside}  to sum over the zeros, gives the weight polynomial as a single residue (or constant term)
\begin{equation}\label{eq:tmtxpower}
(T_L^t)_{y',y} = \res \left [ \frac{P_{y'}(1/x)P^{(y+1)}_{L-y}(1/x)}{x^{t+2}P_{L+1}(1/x)} , \{x, 0\} \right].
\end{equation}
As noted above, a factor, $h_{y',y}$ needs to be inserted in the numerator for the case $y'>y$. Comparing this form with \eqref{eq:classCT} and changing to the reciprocal polynomials \eqref{eq:recippolys} gives  Viennot's theorem. Thus we see the change to the reciprocal $1/x$  in this context   corresponds to switching from a sum of residues to a single residue at infinity.

\subsection*{Using partial fractions}

We can also derive a residue or constant term expression for the weight polynomial, \eqref{eq:tmtxpower},   without  invoking either of the calculus Lemmas \ref{lem:asdasd} or Lemma \ref{lem:insideoutside} by using a partial fraction expansion of a rational function. In particular, if we have the rational function
\[
	G(x)=\frac{Q(x)}{T(x)},
\]
with $Q$ and $T$ polynomials of degree $a$ and $b$ respectively and $a<b$. Assuming $T(x)$ is monic and has simple zeros $T(x_i)=0$  we have $T(x)=\prod_i(x-x_i)$ and thus we have, using standard methods, the partial fraction expansion
\[
	G(x)=  \sum_i \frac{ Q(x_i)}{T'(x_i)}\, \frac{1}{x-x_i}
\]
where $T'(x)$ is the derivative of $T(x)$ and thus
\begin{equation}\label{eq:pfracgf}
		G(1/x)=  \sum_i \frac{ Q(x_i)}{T'(x_i)}\, \frac{x}{1-x_i\,x}.
\end{equation}
Geometric expanding each term gives us the coefficient of $x^n$ in $G(1/x)$ as
\begin{equation}\label{eq:inpfrac} 
	[x^n]G(1/x) = \sum_i x_i^{n-1}\frac{ Q(x_i)}{T'(x_i)}.
\end{equation}
If we now compare \eqref{eq:evexpand} with \eqref{eq:inpfrac} we see that $Q\to P_{y'}  P^{(y+1)}_{L-y}$, $T\to P_{L+1}$ and $t=n-1$, thus we  get 
\begin{align}
	(T_L^t)_{y',y}& =[x^{t+1}] \frac{ P_{y'}(1/x ) P^{(y+1)}_{L-y}(1/x ) }{P_{L+1}(1/x)}\notag\\
	&=\res \left[ \frac{1}{x^{t+2}}   \frac{ P_{y'}(1/x ) P^{(y+1)}_{L-y}(1/x ) }{P_{L+1}(1/x)},\{x,0 \} \right].
\end{align}
Note, the orthogonal polynomials satisfy the conditions required for the existence of \eqref{eq:pfracgf}.
Thus we see that the sum over the zeros (ie.\ eigenvalues) in \eqref{eq:evexpand} is actually a term in the geometric expansion, or Taylor series, of a partial fraction expansion and hence `summed' by reverting the expansion back to the rational function it arose from.  

Note, although the partial fraction route is elementary,   to get  from the natural representation of the $t^\text{th}$ power of a matrix in terms of its diagonalization  (ie.\ equations \eqref{eq:TVDU} and \eqref{eq_quadortho}) to Viennot's theorem we still needed    Lemma \ref{lem:nontrivialPavingPolyIdentity}.

\section{Proof of the Constant Term Theorem}	
\label{sec_resProof}
 
For the proof of the CT Theorem   we use the following residue `change of variable' lemma. 
\begin{lemma}[Residue change of variable] 
\label{lem:changeVars}
 Let $f(x)$  and $r(x)$ be a functions  which have Laurent series about the origin and the Laurent series of $r(x)$ has the property that $r(x)=x^k g(x) $ with $k>0$ and $g(x)$ has a Taylor series  $g(x)=\sum_{n\ge 0}a_n x^n$ such that  $a_0\ne0$. 	Then 
\begin{equation}
 \res \left [ f(x),x\right ] =\frac{1}{k} \res \left [f\bigl(r(z)\bigr)\frac{dr}{dz} , z\right ]
\label{eq_varchange}
\end{equation}
where $\res \left [ f(x), x\right ]$ denotes the residue of $f(x)$ at $x=0$.
\end{lemma}
The lemma appears to have first been proved by Jacobi \cite{jacobi:1830ai}. A proof may also be found in Goulden and Jackson \cite{goulden:1986kx}, Theorem 1.2.2. Note, the condition on $r(x)$ is equivalent to the requirement that the Laurent series of the inverse function of $r(x)$ exist.
	
To prove the CT Theorem  we  start with Viennot's theorem and use   the simple fact that the coefficient of $x^t$ in the Taylor series for $f(x)$ is given by
\begin{equation}\label{eq:ctres}
	 \CT_x\left[\frac{1}{x^t}f(x)\right]=\res\left[\frac{1}{x^{t+1}}f(x)\right]
\end{equation}
which gives \eqref{eq:classCT}, that is,
\begin{equation}\label{eq:classCT22}
	Z_t(y',y;L)= \CT_x\left[\frac{1}{x^{t+1}}x^{y-y'}\,\frac{\recip{P}_{y'}(x)\,h_{y',y}\, \recip{P}^{(y+1)}_{L-y}(x)}{\recip{P}_{L+1}(x)}\right].
\end{equation}
Now consider the change of variable   defined   by
\begin{equation} 
x(\p) = \frac{\p}{\p^2+b\p+ \la },\qquad \lambda\ne0
\label{eq_varsubs1}
\end{equation}
which has the Taylor expansion about the origin
\[
	x(\p) =\frac{\rho }{\lambda }-\frac{b \rho ^2}{\lambda^2}+O\left(\rho ^3\right) 
\] 
and thus $x(\p)$ satisfies the conditions of   Lemma~\ref{lem:changeVars} with  $k=1$. Note, \eqref{eq_varsubs1} is the reciprocal of the change given in \eqref{eq:mup}.  Thus from  \eqref{eq:classCT22} and using   Lemma~\ref{lem:changeVars} with the change of variable \eqref{eq_varsubs1} and derivative
\begin{equation} \label{eq:dxdp}
\frac{d}{d\p}x(\p) = \frac{\la \p^{-2}-1}{(\p+b+\la \p^{-1})^2}.
\end{equation}
we get  the CT theorem result  \eqref{eq:CTrho}, in terms of the Laurent polynomials defined by \eqref{eq:PR}.
 
We make two remarks. The first is the primary reason for the change of variable in this instance is that it ``gets rid of the square roots'' such as those that appear in the representation \eqref{eq_chebEx} and hence changes a representation of a             Chebyshev type polynomial in terms of algebraic functions to an explicit  Laurent polynomial form.

For the second remark we note that for this proof of the CT Theorem we could equally well have started just before  the end of the diagonalization proof of Viennot theorem, ie.\ equation \eqref{eq:tmtxpower}, and proceeded with the residue change of variable \eqref{eq_varsubs1}. In other words, the transfer matrix diagonalization naturally ends with an expression for the weight polynomial -- this has to be generated on before getting to Viennot's theorem,  whilst the first step in this CT proof, ie.\ equation \eqref{eq:classCT22},  is to undo this generating step.


%
\section{Proof of the Paving Theorem}	
%

\begin{figure}[htbp]
	\begin{center}
	\includegraphics[width=10cm, height=20mm]{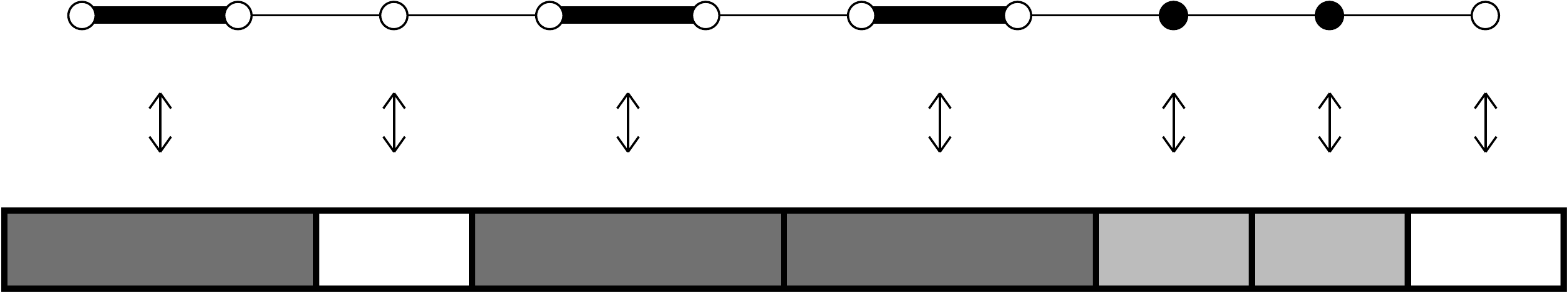}
	\end{center}
	\caption{The upper part of this diagram shows a \defin{paving}.  The lower part indicates the reason for calling it a `paving', as it is in bijection with a more standard `paving' or `tiling' diagram, where long and short tiles are used, the short tiles being of two possible colours.}
	\label{fig:pavingEx2}

\end{figure}

The paving interpretation of the orthogonal polynomial three term recurrence relation was introduced by Viennot \cite{Viennot1985ah}. We will use this interpretation as the primary means  of proving Theorem~\ref{thm:ViennotPaving}. First we define several terms associated with pavings.  A \defin{path graph} is any graph isomorphic to a graph with vertex set $\{v_i\}_{i=0}^k$ and edge set $\{v_iv_{i+1}\}_{i=0}^{k-1}$.  
A \defin{monomer} is a distinguished vertex in a graph.  A \defin{dimer} is a distinguished edge (pair of adjacent vertices).  A \defin{non-covered vertex} is a vertex which occurs in neither a monomer nor a dimer.  A \defin{paver} is any of the three possibilities:  monomer, dimer or non-covered vertex.  A \defin{paving} is a collection of pavers on a path graph such that no two pavers share a vertex.  
We say that a paving is \defin{order $k$} if it occurs on a path graph  with $k$ vertices.
An example of a paving of order ten is shown in Figure~\ref{fig:pavingEx2}.  

\defin{Weighted pavings} are pavings with weights associated with each paver.
We will need pavers with shifted indices in order to calculate the shifted paving polynomials that occur in Theorem~\ref{thm:ViennotPaving}.
Thus, the \defin{weight of a paver} $\alpha$ with \defin{shift} $j$ is defined as follows. 
\begin{equation} \label{eq:weightPaver}
w_j(\alpha)=
\begin{cases}
x & \text{the paver $\alpha$ is a non-covered vertex}\\
-b_{i+j} & \text{the paver $\alpha$ is the monomer $v_i$}\\
-\la_{i+j} & \text{the paver $\alpha$ is the dimer $e_i$}
\end{cases}
\end{equation}
The \defin{weight of a paving} is defined to be the product of the weights of the pavers that comprise it, i.e.
\begin{equation}
w_j(p) = \prod_{\alpha \in p}w_j(\alpha),
\end{equation}
for $p$ a paving.  It is useful to distinguish two kinds of paving.
Pavings containing only non-covered vertices and dimers are called \defin{Ballot pavings}; those also containing monomers are called \defin{Motzkin pavings}.  A \defin{paving set} is the collection of all pavings (of either Ballot or Motzkin type) on a path graph of given size.  We write
\begin{align}
	\mathcal{P}^{\text{Bal}}_{k} &=  \{p | \text{$p$ is a Ballot Paving of order $k$}\}
	\label{eq:PavingSetBal}\\
	\mathcal{P}^{\text{Motz}}_{k} &=  \{p | \text{$p$ is a Motzkin Paving of order $k$}\}
	\label{eq:PavingSetMotz}
\end{align}
When it is clear by context whether we refer to sets of Ballot or Motzkin pavings, the explanatory superscript is omitted.

A \defin{paving polynomial} is a sum over weighted pavings defined by
\begin{equation} \label{eq:pavingPoly}
P_k^{(j)}(x) = \sum_{p \in \mathcal{P}_k} w_j(p),
\end{equation}
with $\mathcal{P}_k$   being either $\mathcal{P}^{\text{Bal}}_k$ or $\mathcal{P}^{\text{Motz}}_k$ and weights as in Equation~(\ref{eq:weightPaver}).  If $k=0$ then we define $P_0^{(j)}(x)\df 1$ $\forall j$. 
The diagrammatic representation for a paving polynomial on a  paving set of Ballot type is
\begin{equation} \label{eq:setdiagBal_5}
\includegraphics[width=11cm]{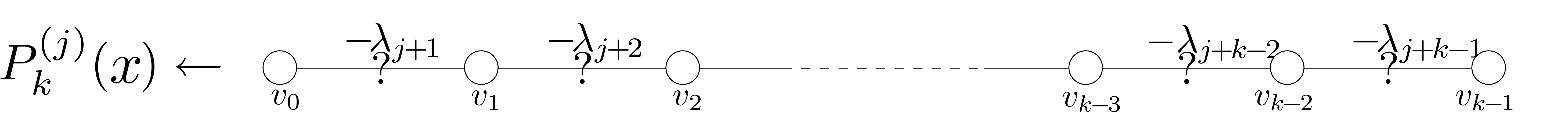}
\end{equation}
where the question mark denotes that the edge can be either a dimer or not.
The diagrammatic representation for a paving polynomial on a   paving set of Motzkin type is 
\begin{equation} \label{eq:setdiagMotz_5}
\includegraphics[width=11cm]{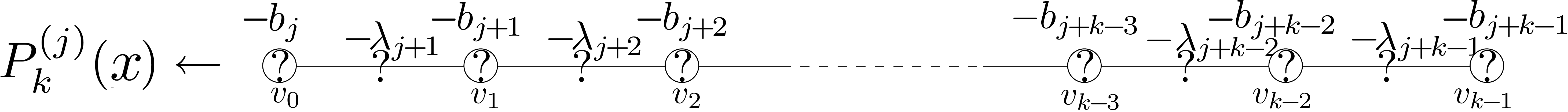}
\end{equation}
Note, we overload the notation for $P_k^{(j)}(x)$ since we will use it for the set of pavings and the corresponding paving polynomial obtained by summing over all weighted pavings in the set.

Viennot \cite{Viennot1985ah}  has shown that Equation~(\ref{eq:pavingPoly}) satisfies  \eqref{eq:threetermrec}  and hence $P_k^{(j)}(x)$ is  an orthogonal  polynomial.  Ballot Pavings correspond to the case $b_i=0$ $\forall i$.
An example of   Ballot pavings and the resulting polynomial  is
\begin{center}
\includegraphics[width=10cm]{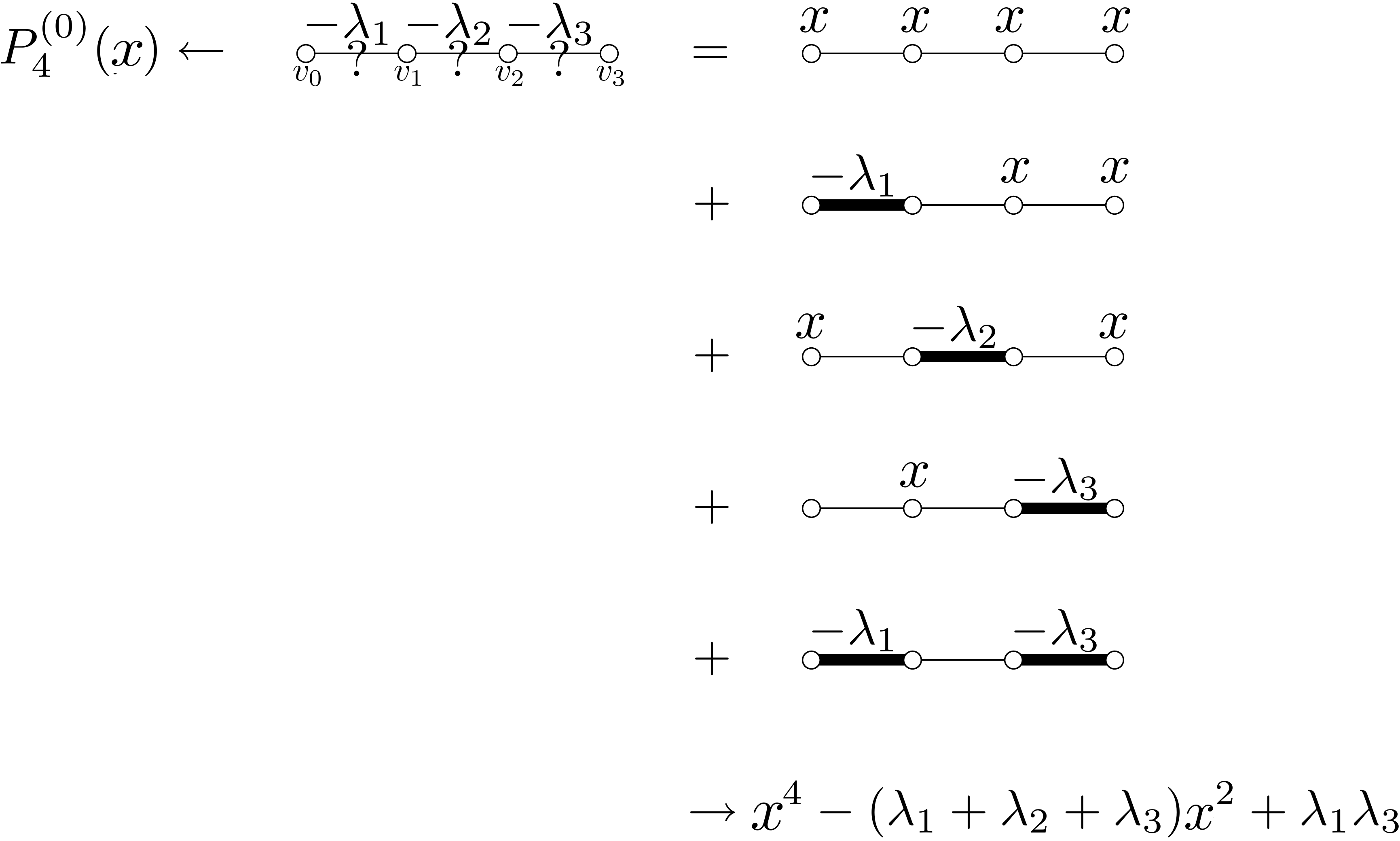}
\end{center}
and an example of a Motzkin paving set with associated polynomial  is
\begin{center}
\includegraphics[width=11cm]{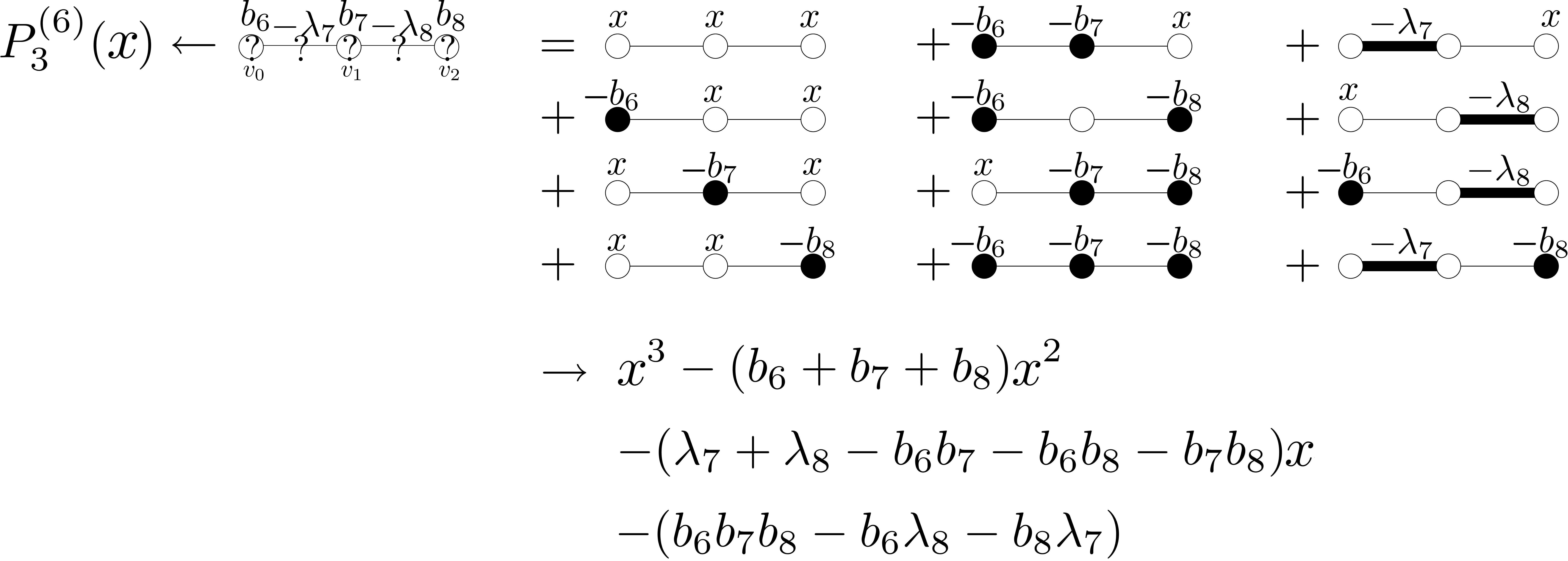}
\end{center}

The two identities of Theorem~2 correspond to
\begin{itemize}
	\item cutting at an arbitrary edge; and 
	\item cutting at an arbitrary vertex.
\end{itemize}
In the first procedure we consider the edge $e_c$.  This edge is either paved or not paved with a dimer.   Equation~(\ref{eq:setdiagMotzCutAtEdgeC}) illustrates the division into these two cases, and the corresponding polynomial identity.  Note that the right hand side of the expression obtained is a sum of products of smaller order polynomials such that the weight `$-\lambda_{c+j}$' associated with edge $e_c$ occurs explicitly as a coefficient; and is not hidden inside any of the smaller order polynomials.
\begin{equation} \label{eq:setdiagMotzCutAtEdgeC}
\includegraphics[width=10.5cm]{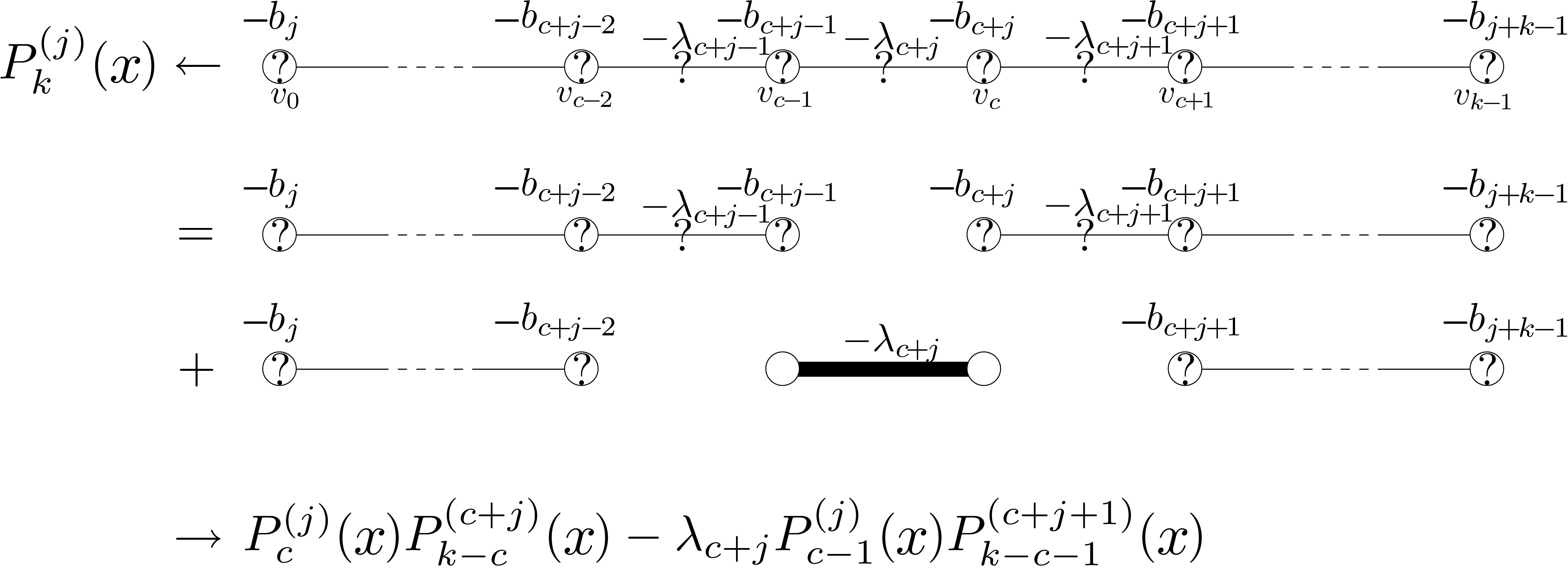}
\end{equation}
The second procedure is to cut at an arbitrary vertex `$v_c$'.  In this procedure the cases to consider are: $v_c$ is non-covered, $v_c$ is a monomer, $v_c$ is the leftmost vertex of a dimer, and $v_c$ is the rightmost vertex of a monomer.  These four cases are shown in Equation~(\ref{eq:setdiagMotzCutAtVertexC}), and the resulting identity gives $P_k^{(j)}(x)$ as a sum of four terms, each of which contains a product of two smaller order polynomials.  The vertex weight 
`$-b_{c+j}$' occurs explicitly as a coefficient, and is not hidden in any of the smaller order polynomials.
\begin{equation} \label{eq:setdiagMotzCutAtVertexC}
\includegraphics[width=11cm]{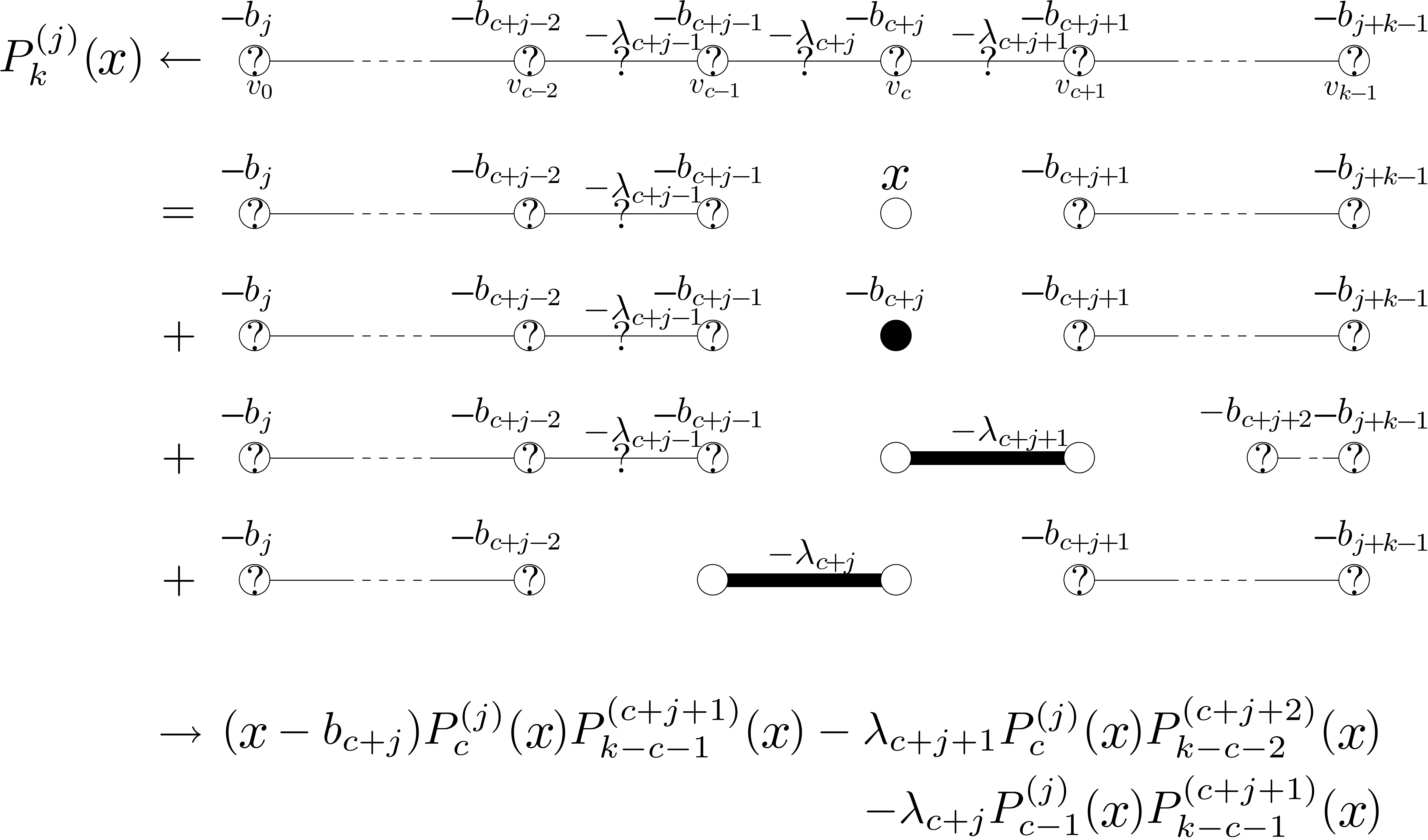}
\end{equation}
These two `cutting' procedures, as shown in Diagrams~(\ref{eq:setdiagMotzCutAtEdgeC}) and (\ref{eq:setdiagMotzCutAtVertexC}), prove Part~\ref{it:oneCut} of Theorem~\ref{thm:ViennotPaving}.
Finally Part~\ref{it:bound} of the Theorem follows immediately from Part~\ref{it:oneCut} by induction on the number of decorations of each of the two possible kinds (`across step' and `down step').

We note that the proof of Theorem~\ref{thm:ViennotPaving} shows that the upper bounds given in Equations~(\ref{eq:bounds}) are tight only when the decorations are well-separated from each other as well as from the ends of the diagram.  Otherwise, pulling out a given decoration as a coefficient may pull out a neighboring decoration in the same procedure, meaning that fewer terms are needed.  Thus we can deal more efficiently with decorated weightings in which collections of decorations are bunched together, as illustrated in the examples in Section~\ref{sec:applications}.

\section{Applications}\label{sec:applications}  
\label{sec:apps}

We now consider two applications. The first is the DiMazio and Rubin problem discussed in the introduction. Solving this model  corresponds to determining the Ballot path weight polynomials with just one upper and one lower decorated weight. The second application is an extension of that problem in which two upper and two lower edges now carry decorated weights.

\subsection{Two decorated weights}	

For the DiMazio and Rubin problem we need to compute the Ballot path weight polynomials with weights given by \eqref{eq_dirub}. A partial solution was given in  \cite{brak:2005fa} for the case when the upper weight is a particular function of the lower one ie.\ $\kappa+\omega=\kappa\omega$.  The first general solution was published in \cite{brak:2006kx}  in 2006.
\begin{figure}[htbp] 
\begin{center}
\includegraphics[width=10cm]{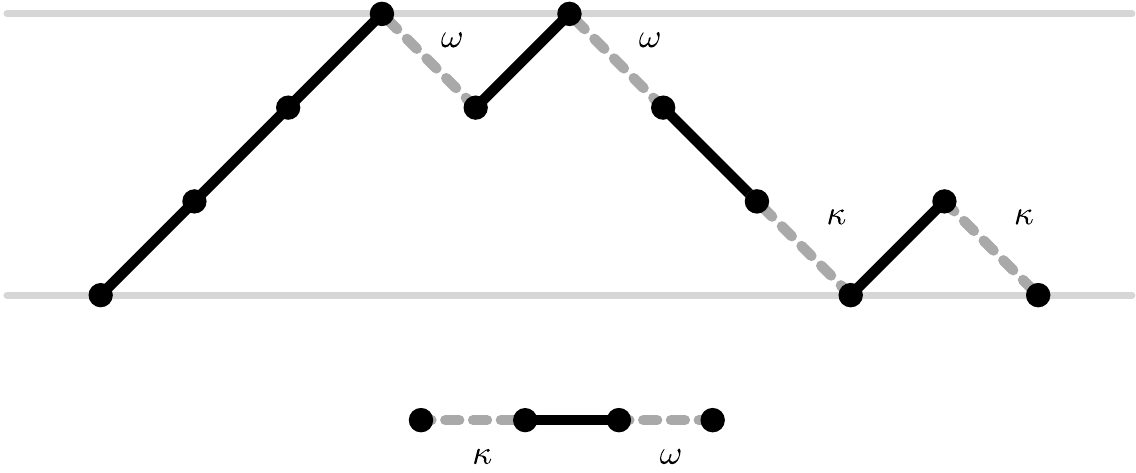}
\caption{An example of a two weight Dyck path (above) in a strip of height three and the corresponding paving problem (below).}
\label{fig:dyckOneWeight}
\end{center}
\end{figure}
The solution we now give is an improvement on  \cite{brak:2006kx}  (which was based on a precursor to this method) as it contains fewer summations.  

\begin{theorem} \label{thm:DMCT}  
	Let $Z_{2r}(\ka, \w;L)$ be the weight polynomial for Dyck paths of length $2r$ confined to a strip of height $L$, and with weights (see \eqref{eq:lailall}),  $b=\lambda=1$, $\hat b_i=0\ \forall i$ and
\begin{align}
	\hat\lambda_i=\begin{cases}
		\omega-1 & \text{if $i=L$}\\
		\kappa-1 & \text{if $i=1$}\\
		0 & \text{otherwise.}\\
	\end{cases}
\label{eq_dirub}
\end{align}		
%
Then 
\begin{equation}
Z_{2r}(\kappa,\omega;L) =\CT \left [ (\p+\p^{-1})^{2r} (1-\p^2) \frac{A\p^L-B\p^{-L}}{AC\p^L -BD \p^{-L}}\right ]
\label{eq_thmtwoweights}
\end{equation}
where
\begin{subequations} \label{eq:ABCD}
\begin{align}
A &=  \p^2 - \hat{\omega} \\
B &=  1 -\hat{\omega} \p^2 \\
C &=  \p^2 - \hat{\kappa} \\
D &=  1 -\hat{\kappa} \p^2, 
\end{align}	
\end{subequations}
$\hat{\ka} \df  \ka -1$ and $\hat{\w} \df  \w - 1$.
\end{theorem}
An example of the paths referred to in Theorem~\ref{thm:DMCT} is given in Figure~\ref{fig:dyckOneWeight}. 
Theorem~\ref{thm:DMCT} is derived using Theorem~\ref{thm:CT} and Theorem~\ref{thm:ViennotPaving}.   We shall show the derivation of  the denominator of \eqref{eq_thmtwoweights} but omit the details for  the numerator which may be similarly derived.  

We work with polynomials in the $x$ variable, and then use Equation~(\ref{eq:PR}) to get the expression in terms of $\p$. First write $P_{L+1}(x)$ for our given weighting as in Equation~(\ref{eq:dmPoly1new}):
\begin{equation} \label{eq:dmPoly1new}
	\includegraphics[width=8.5cm]{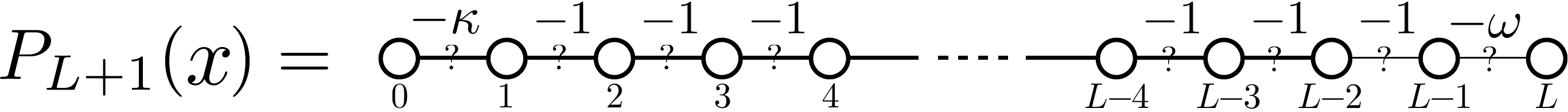}	
\end{equation}
We cut at the first and last edges, as in Equation~(\ref{eq:dmPoly2new}).  (The edges to cut at are chosen since they mark the boundary between decorated and undecorated sections of the path graph.)
\begin{equation} \label{eq:dmPoly2new}
\includegraphics[width=8.5cm]{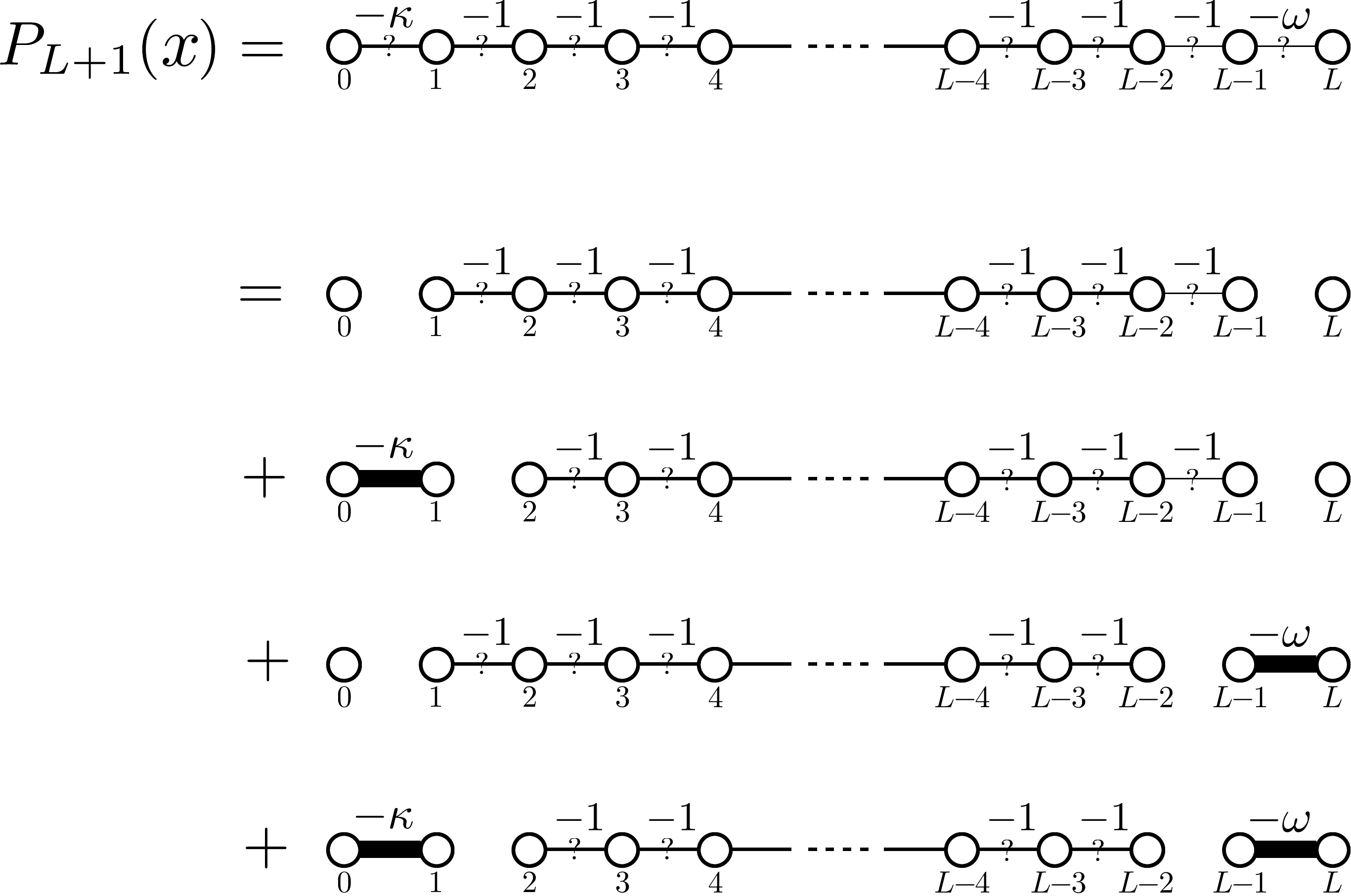} 
\end{equation}
Now we have $P_{L+1}(x)$ as a sum of four terms, as in Equation~(\ref{eq:dmPoly3new}):
\begin{equation} \label{eq:dmPoly3new}
P_{L+1}(x) = x^2 S_{L-1}(x) - \kappa x S_{L-2}(x) - \omega x S_{L-2}(x) + \kappa \omega S_{L-3}(x).
\end{equation}
Each of the four terms is a product of three contributions - the $x$'s, $\kappa$'s and $\omega$'s come from the short sections in each row of the diagram, and $S_k$'s represent the long sections, as in Equation~(\ref{eq:sPolynew}):
\begin{equation} \label{eq:sPolynew}
\includegraphics[width=7.5cm]{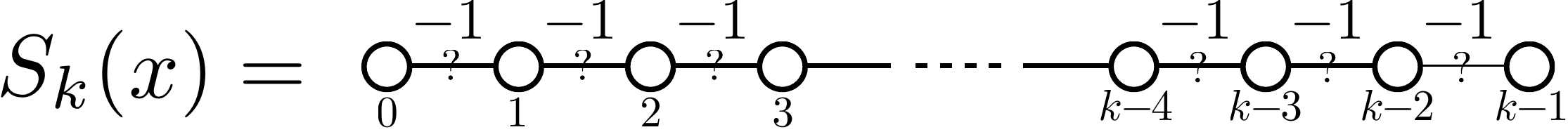}
\end{equation}
Now Equation~(\ref{eq:sPolynew}) represents polynomials satisfying the recurrence and initial conditions given in Equation~(\ref{eq:Sk}), so we may substitute Equation~(\ref{eq_chebEx}) for each occurrence on an `$S_k$' in Equation~(\ref{eq:dmPoly3new}).  We now have a sum of ratios of surds, which may be simplified by the change of variables specified in Equation~(\ref{eq:PR}) to yield
\begin{align}
R_{L+1}(\p) &= P_{L+1}(\p+  \p^{-1}) \\
&=(AC\p^L - BD \p^{-L})\p^2(\p-\p^{-1})
\end{align}
which is, up to a factor, the denominator of \eqref{eq_thmtwoweights}. The numerator is similarly derived.
 
Next we indicate how to expand Equation~(\ref{eq_thmtwoweights}) to give an expansion for the weight polynomial in terms of binomials.  Our particular expression  is obtained via straightforward application of geometric series and binomial expansions, with the resulting formula containing a 5-fold sum.  It is certainly possible to do worse than this and obtain more sums by making less judicious choices of representations while carrying out the expansion, but it seems unlikely that elementary methods can yield a smaller than 5-fold sum for this problem.

First, manipulate the fractional part of Equation~(\ref{eq_thmtwoweights}) into a form in which geometric expansion of the denominator is natural, as in Line~(\ref{eq:geomReady}) below.  Do the geometric expansion and multiply out to give two terms: extract the $m=0$ case from the first term and shift the index of summation in the second to give Line~(\ref{eq:threeParts}).
\begin{align}
\frac{A \p^L - B \p^{-L}}{AC \p^L - BD \p ^{-L}}
&=  
\frac{\frac{1}{D}  - \frac{A}{BD} \p^{2L}}{1  - \frac{AC}{BD} \p^{2L}} \label{eq:geomReady}\\
\nonumber \\
&= 
\left(\frac{1}{D}  - \frac{A}{BD} \p^{2L} \right )
\sum_{m=0}^{\infty} \left ( \frac{AC}{BD} \right )^m\p^{2mL}\\
\nonumber \\
&=  \frac{1}{D} + 
\sum_{m=1}^{\infty}  \frac{A^mC^m}{B^mD^{m+1}} \p^{2mL} -
\sum_{m=1}^{\infty}  \frac{A^{m}C^{m-1}}{B^{m}D^{m}} \p^{2mL} 
\nonumber \label{eq:threeParts}\\
\end{align}
We next find the $\CT$, separately, of each of the three terms multiplied by series for $(\p+\p^{-1})^{2r}(1-\p^2)$.  The first gives

\begin{align}
  \text{CT} &\left [ (\p+\p^{-1})^{2r} (1-\p^2)
\frac{1}{D}\right ]\\
\nonumber \\
&=  \text{CT} \left [\left (\sum_{u=0}^{2r}\binom{2r}{u}\p^{2r-2u}\right ) (1-\p^2)
\left (\sum_{m=0}^{\infty} \hat{\ka}^m \p^{2m}\right )\right ]\\
\nonumber \\
&=  \sum_{m=0}^{\infty} \hat{\ka}^m \text{CT} \left [\sum_{u=0}^{2r}\binom{2r}{u}(\p^{2r-2u+2m}-\p^{2r-2u+2m+2})\right ]\\
\nonumber \\
&=  \sum_{m=0}^{\infty} \hat{\ka}^m \left [\binom{2r}{r+m}-\binom{2r}{r+m+1}\right ]
\\
\nonumber \\
&=  \sum_{m=0}^{\infty}C_{r;r-m} \hat{\ka}^m 
\end{align}
where $C_{n;k}:=	\binom{2n}{k}-\binom{2n}{k-1}$.
The second term is expanded similarly,
with the positive powers of $A$ and $C$ and the negative powers of $B$ and $D$ each contributing a single sum, which, when concatenated with the original sum over $m$, creates a 5-fold sum altogether.  The third term generates a similar 5-fold sum; and this difference of a pair of 5-fold sums is combined into one in the final expression in Theorem~\ref{thm:DMCTexpanded} below.  

\begin{theorem} \label{thm:DMCTexpanded} 
	Let $Z_{2r}(\kappa,\omega;L)$ be as in Theorem~\ref{thm:DMCT}.  Then
\begin{equation} 
\begin{split}
Z_{2r}(\kappa,\omega;L) =& \sum_{m \ge 0} C_{r; r-m}\hat{\ka}^m  \\
&
+ \sum_{m \ge 1} \sum_{p_1, p_2 \ge 0} \sum_{s_1, s_2 =0}^m (-1)^{s_1+s_2} \hat{\ka}^{s_2+p_2} \hat{\w}^{s_1 + p_1}\\
& \times
\binom{m}{s_1}\binom{m}{s_2} \binom{m-1+p_1}{p_1}\binom{m+p_2}{p_2}\\
& \times \bigg [
C_{r;r-k-1}   
- \frac{m-s_2}{m+p_2}C_{r;r-k}  
\bigg ]
\label{eq_tw}
\end{split}
\end{equation}
where $k=p_1+p_2-s_1-s_2+(L+2)m-1$, and $C_{n;k}$ is the extended Catalan number,
\[
	C_{n;k}=	\binom{2n}{k}-\binom{2n}{k-1}.
\]
The binomial coefficient $\binom{n}{m}$ is assumed to vanish if $n< 0$ or  $m <0 $ or $n<m$.
\end{theorem}
Note, when trying to rearrange \eqref{eq_tw} care should be taken when using any binomial identities   because of the vanishing condition on the binomial coefficients - the support of any new expression must be the same as the support before (alternatively the upper limits of  all of the summations must be precisely stated).


\subsection{Four Decorated weights}	

\begin{figure}[htbp]
\begin{center}
\includegraphics[width=10cm]{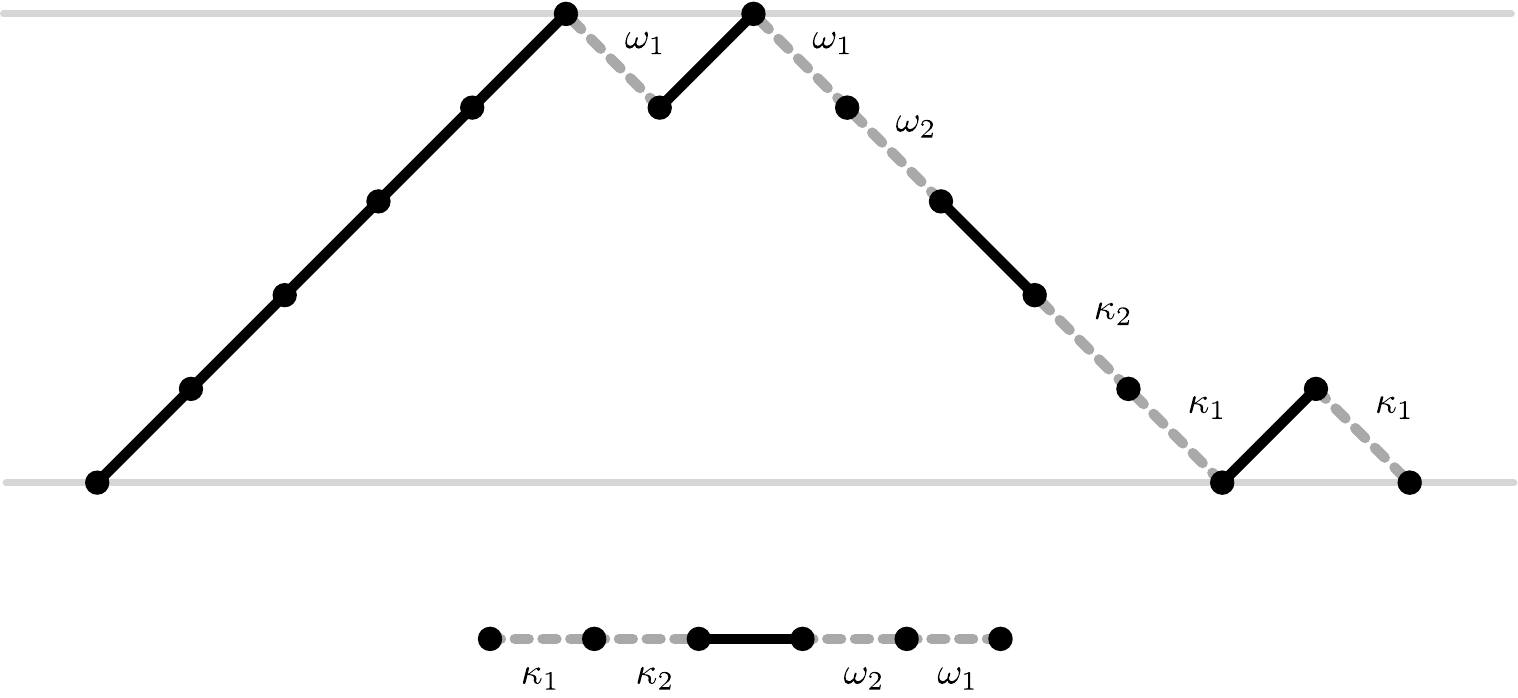}
\caption{An example of a four weight Dyck path (above) in a strip of height five and the corresponding paving problem (below).}
\label{fig:dyckTwoWeights}
\end{center}
\end{figure}

This second problem is a natural generalization of the previous problem.  We now have a pair of decorated  weights in the pair of rows adjacent to each wall, as in Figure~\ref{fig:dyckTwoWeights}.
In the earlier DiMazio-Rubin problem, paths have been interpreted as polymers zig-zagging between comparatively large colloidal particles (large enough to be approximated by flat walls above and below) with an interaction occurring only upon contact between the surface and the polymer; this weighting scheme could be used to model such polymer systems, but now with a longer range interaction strength that varies sharply with separation from the colloid.
\begin{theorem} \label{thm:2pairweightsCT}  Let $Z_{2r}(\ka_1, \ka_2, \w_1, \w_2;L)$ be the weight polynomial for Dyck paths of length $2r$ confined to a strip of height $L$,   
and with weights (see \eqref{eq:lailall}),  $b=\lambda=1$, $\hat b_i=0\ \forall i$ and
\begin{equation}
	\hat\lambda_i=\begin{cases}
		\omega_1-1 & \text{if $i=L$}\\
		\omega_2-1 & \text{if $i=L-1$}\\
		\kappa_2-1 & \text{if $i=2$}\\
		\kappa_1-1 & \text{if $i=1$}\\
		0 & \text{otherwise.}
	\end{cases}
\label{eq_dirubext}
\end{equation}
Then 
\begin{align} \label{eq:CT2pairweights}
&Z_{2r}(\ka_1, \ka_2, \w_1, \w_2;L)\notag \\
&=\CT \left [ (\p+\p^{-1})^{2r} \left ( \frac{AB \p^L - \Abar \,\Bbar \p^{-L}}{C B \p^L - \Cbar\, \Bbar \p^{-L}}\right ) (\p^{-1}-\p)\right ]
\end{align}
where
\begin{subequations}
\begin{align}
A &= 1 - \hat{\ka}_2 \p^{-2}\\
\Abar &= 1 - \hat{\ka}_2 \p^{2}\\
B &= \p - (\hat{\w}_1+\hat{\w}_2)\p^{-1} -\hat{\w}_2\p^{-3}\\
\Bbar &= \p^{-1} - (\hat{\w}_1+\hat{\w}_2)\p -\hat{\w}_2\p^3\\
C &= \p - (\hat{\ka}_1+\hat{\ka}_2)\p^{-1} -\hat{\ka}_2\p^{-3}\\
\Cbar &= \p^{-1} - (\hat{\ka}_1+\hat{\ka}_2)\p -\hat{\ka}_2\p^{3}
\end{align}	
\end{subequations}
for $\hat{\ka}_i \df  \ka_i -1$ and $\hat{\w}_i \df  \w_i - 1$.
\end{theorem}
An example of the paths referred to in the above theorem is given in Figure~\ref{fig:dyckTwoWeights}.

The constant term expression in Theorem~\ref{thm:2pairweightsCT} may be expanded in a similar fashion to that of Theorem~\ref{thm:DMCT} to yield a 9-fold sum.  The fractional component of Equation~(\ref{eq:CT2pairweights}) may be written 

\begin{align}
\frac{AB\p^L-\Ao \Bo \p^{-L}}{C B \p^L - \Co \Bo \p^{-L}}
&=& \frac{\Ao}{\Co} +
\Ao \sum_{m=1}^\infty \frac{C^m B^m}{\Co^{m+1} \Bo^m}\p^{2mL} 
-
A \sum_{m=1}^\infty \frac{C^{m-1} B^{m}}{\Co^{m} \Bo^{m}} \p^{2mL} \nonumber \label{eq:4weightsfractional}\\
\end{align}
by a method precisely analogous to that applied in the 2-weights case.  When multiplied by $(\p+\p^{-1})^{2r}(\p^{-1}-\p)$ and the constant term extracted, the initial term gives the double summation in the first term of Equation~(\ref{eq:4weightsExpansion}) in Theorem~\ref{thm:twoWeightsPerWallExpansion} below.  The other two terms of Equation~(\ref{eq:4weightsfractional}) each give 9-fold sums, as a consequence of the double sums yielded by each of the powers, positive and negative, of $C$, $B$, $\Co$ and $\Bo$.  These two 9-fold sums are combined in the second term of Equation~(\ref{eq:4weightsExpansion}) in Theorem~\ref{thm:twoWeightsPerWallExpansion}.
\begin{theorem} \label{thm:twoWeightsPerWallExpansion} Let $Z_{2r}(\ka_1, \ka_2, \w_1, \w_2;L)$ be as in Theorem~\ref{thm:2pairweightsCT}.  Then
\begin{small}
\begin{align}
Z_{2r}&(\ka_1, \ka_2, \w_1, \w_2;L)\notag \\
=&\sum_{i\ge0} \sum_{j=0}^i \binom{i}{j} (\hat{\ka}_1 + \hat{\ka}_2)^j \hat{\ka}_2^{i-j} \times \notag \\
&\left ( \hat{\ka}_2 \binom{2r}{u_0+2} -(\hat{\ka}_2+1) \binom{2r}{u_0+1} + \binom{2r}{u_0}\right )\notag \\
+&\sum_{m \ge 1} 
\sum_{s_1=0}^m \sum_{i_1=0}^{s_1} \sum_{s_2=0}^m \sum_{i_2=0}^{s_2} \sum_{v_1 \ge 0} \sum_{j_1=0}^{v_1} \sum_{v_2 \ge 0} \sum_{j_2=0}^{v_2}
\binom{s_1}{i_1}\binom{m}{s_2}\binom{s_2}{i_2} 
\binom{v_1}{j_1} \times \notag \\
&\binom{v_2+m-1}{m-1}\binom{v_2}{j_2}(-1)^{s_1+s_2+i_1+i_2} \times\notag \\
&\hat{\ka}_2^{i_1+j_1}\left ( \hat{\ka}_1 + \hat{\ka}_2 \right )^{m+v_1-1-s_1-j_1}
\hat{\w}_2^{i_2+j_2}\left ( \hat{\w}_1 + \hat{\w}_2 \right )^{m+v_2-s_2-j_2}\times\notag \\
&{\bigg \{}
\binom{m}{s_1} \binom{v_1+m}{m}(\hat{\ka}_1 + \hat{\ka}_2)
\left ( \hat{\ka}_2\binom{2r}{u_1+2}-(\hat{\ka}_2+1)\binom{2r}{u_1+1}+\binom{2r}{u_1}\right )\notag \\
&-\binom{m-1}{s_1} \binom{v_1+m-1}{m-1} 
\left ( \hat{\ka}_2\binom{2r}{u_1-1}-{(\hat{\ka}_2+1)}\binom{2r}{u_1}+\binom{2r}{u_1+1}\right )
{\bigg \}} \label{eq:4weightsExpansion}
\end{align}
\end{small}

for $u_0=r+2i-j$ and $u_1 = r+mL+v_1+v_2+s_1+s_2+j_1+j_2-2i_1-2i_2$.
\end{theorem}

Theorems~\ref{thm:DMCTexpanded} and \ref{thm:twoWeightsPerWallExpansion} are to be compared with the Rogers formula below, which gives the weight polynomial as an order $L$-fold sum.

\begin{theorem}[Rogers \cite{rogers:1907rc}] \label{thm:allweightsDyck} Let $Z_{2n}$ be the weight polynomial for the set of Dyck paths of length $2n$ with general down step weighting (see \eqref{eq:lailall}):  $b=\lambda=1$, $\hat b_i=0\ \forall i$ and $\hat\lambda_i= \kappa_i-1$	
in either a strip of height $L$ or in the half plane (take $L=\infty$).  Then the weight polynomial is given by
\begin{equation}
Z_{2n}(\ka_1, \ka_2, ...;L)=\sum_{l=0}^{\min \{n-1, L-1\}} s_l,
\end{equation}
where $s_l$ is the weight polynomial for that subset of paths in $Z_{2n}$ which reach but do not exceed height $L+1$.   The $s_l$'s are given by
\begin{equation}
s_0=\ka_1^n,
\end{equation}
and
\begin{equation}
s_l = 
\sum_{j_1=l}^{j_0-1}\sum_{j_2=l-1}^{j_1-1} \hdots \sum_{j_l=1}^{j_{l-1}-1} 
\prod_{k=0}^{l-1} \binom{j_k-j_{k+2}-1}{j_{k}-j_{k+1}} \ka_1^{j_0-j_1} \ka_2^{j_1-j_2} \hdots \ka_{l+1}^{j_l-j_{l+1}},
\end{equation}
for $l \ge 1$; with
\begin{align}
j_0 &\df   n, \\
j_{l+1}&\df   0.
\end{align}
\end{theorem}

We see that there is a trade-off between having comparatively few sums (compared with the width of the strip) but a complicated summand, as in Theorems~\ref{thm:DMCTexpanded} and \ref{thm:twoWeightsPerWallExpansion}, versus having a simpler summand but the order of $L$ sums irrespective of the number of decorations.

It would be an interesting piece of further research to see whether the solutions to the problems presented in Section~\ref{sec:applications}, containing as they do multiple alternating sums, are in some appropriate sense best possible or not.  Another potentially useful area of further research would be an investigation of good techniques for extracting asymptotic information directly from the $\CT_\p$ expression.  It would also be interesting to have a pure algebraic formulation of this constant term method -- one that does not rely on any residue theorems.



\section*{Acknowledgements}

Financial support from the Australian Research Council is gratefully
acknowledged. One of the authors, JO thanks the
Graduate School of The University of Melbourne for an Australian
Postgraduate Award and The Centre of Excellence for the Mathematics and Statistics of Complex Systems (MASCOS) for additional financial support. We would also very much like to thank  Ira Gessel for some useful insights and in particular to drawing our attention to the Jacobi and Goulden and Jackson reference \cite{goulden:1983yg} used in the proof in Section \ref{sec_resProof}.    

\bibliographystyle{plain}
\bibliography{CTBibtekDatabase}

\end{document}